\documentclass[12pt]{article}
\usepackage{tikz,float,hyperref,collref}
\usepackage[utf8]{inputenc}
\usepackage[T1]{fontenc}
\usepackage[margin=2.75cm]{geometry}
\usepackage{amsmath,amsfonts,mathtools,authblk,amssymb,amsthm}
\usepackage{cleveref,graphicx,tabularx,ragged2e} 
\usepackage{booktabs,dirtytalk,multicol}
\newtheorem{theorem}{Theorem}[section]     
\newtheorem{corollary}[theorem]{Corollary}

\newtheorem{definition}[theorem]{Definition}
\newtheorem{example}[theorem]{Example}

\allowdisplaybreaks
\date{}

\title{On Laplacian and Distance Laplacian Spectra of Generalized Fan Graph \& a New Graph Class}

\author{\noindent\large Subarsha Banerjee$^{1}$\footnote{Corresponding author.\\ Email address:
		\href{mailto:subarsha.banerjee@jisuniversity.ac.in}{subarsha.banerjee@jisuniversity.ac.in/subarshabnrj@gmail.com}},
	and Soumya Ganguly$^{2}$ }

\affil{$^{1}$\small \footnotesize Department of Mathematics, JIS University, Kolkata, West Bengal 700109, India. \\
	$^{2}$\small \footnotesize  BTech(2nd Year), Department of Computer Science \& Engineering, JIS University, Kolkata, West Bengal 700109, India.}

\begin{document}
	\maketitle

	\begin{abstract}
		Given a graph $G$, the Laplacian matrix of $G$, $L(G)$ is the difference of the adjacency matrix $A(G)$ and $\text{Deg}(G)$, where $\text{Deg}(G)$ is the diagonal matrix of vertex degrees.
		The distance  Laplacian matrix $D^L({G})$ is the difference of the  transmission matrix of $G$ and the distance matrix of $G$.
		In the given paper, we first obtain the Laplacian and distance Laplacian spectrum of generalized fan graphs.
		We then introduce a new graph class which is denoted by $\mathcal{NC}(F_{m,n})$. Finally, we determine the Laplacian spectrum and the distance Laplacian spectrum of $\mathcal{NC}(F_{m,n})$.
		
	\end{abstract}
	
	\textbf{Keywords:}  Laplacian spectrum; distance Laplacian spectrum; generalized fan graph; equitable partition. 
	\\
	\textbf{2010 Mathematics Subject Classification:}  05C07, 05C12, 05C50.
	
	\section{Introduction}
	Throughout the paper, $G$ shall denote a  finite, simple, and undirected graph.
	Let $V(G)=\{v_1,v_2,\dots, v_n\}$ denote the set of all vertices of  $G$, and let $E(G)$ denote the set of all edges of $G$.
	The \textit{order} of $G$ is the number of elements in $V(G)$. 
	Let  $v_i,v_j\in V(G)$. We say that the vertex $v_i$ to be \textit{adjacent} to $v_j$ provided there is an edge from $v_i$ to $v_j$ or vice versa.
	If the vertices $v_i$ and $v_j$ are adjacent to each other, it shall be  denoted  by $v_i\sim v_j$.
	The total number of vertices in $G$ that are adjacent to a given vertex  $v$ is known as the \textit{degree} of $v$.
	The \textit{join} of two graphs $G_1$ and $G_2$ is is denoted by $G_1+G_2$.
	
	The \textit{adjacency} matrix $A(G)$ of $G$ is defined as $A(G)=(a_{ij})_{n\times n}$ is an $n\times n$ matrix defined as follows:
	$a_{ij}=\begin{cases}
	1 & \text{ if } v_i\sim v_j\\
	0 & \text{ elsewhere }.
	\end{cases}$.
	The \textit{Laplacian} matrix $L(G)$  of $G$ is defined as $L(G)=(l_{ij})_{n\times n}$ is defined as follows:
	$l_{ij}=\begin{cases}
	d_i & \textbf{ if } i=j\\
	-1 & \text{ if } v_i\sim v_j\\
	0 & \text{ elsewhere }.
	\end{cases}$.
	Here, $d_i$ denotes the degree of the $i^{th}$ vertex $v_i$.
	The Laplacian matrix $L(G)$ of a graph $G$ has all its eigenvalues as real numbers.
	Moreover, $L(G)$ is a positive semidefinite matrix.
	Consequently, all the real eigenvalues of $L(G)$ are non-negative. 
	It is known that the summation of row  entries in a Laplacian matrix is zero.
	Thus, the determinant of $L(G)$ is always $0$.
	Hence, $0$ is always an eigenvalue of $L(G)$.
	
	A sequence of vertices and edges in a graph $G$ is known as a \textit{walk}. A walk is said to be \textit{closed} if the starting vertex is the same as the end vertex.    
	If all the edges are different in a walk, then it is known as a \textit{trail.} A \textit{path} is a trail in which no vertex is repeated. A closed path is said to be a \textit{cycle}.
	The number of edges in a path is known as the \textit{length} of the path.
	The \textit{distance} matrix of a connected graph $G$ is defined as $D(G)=(d_{ij})_{n\times n}$, where $d_{ij}=d(v_i,v_j)$ is the distance between two vertices $v_i$ and $v_j$.
	The sum of distances from a vertex $v$ to all other vertices of ${G}$ is known as the \textit{transmission} of $v$.
	The transmission of a vertex $v$ is denoted by $Tr(v).$
	The \textit{transmission matrix}  of $G$ is an $n\times n$ matrix where each diagonal entry denotes the transmission of the vertex $v$, and each off-diagonal entry is $0$. 
	The \textit{distance  Laplacian} matrix  $D^L({G})$ of a connected graph  $G$ is defined as  $D^L({G})=Tr({G})-D({G})$. It was introduced in \cite{1}.
	The \textit{distance signless Laplacian} matrix $D^Q({G})$ is defined as  $D^{Q}({G})=Tr({G})+D({G})$.
	Recently, the researchers have studied the two matrices extensively, see for example \cite{2}, \cite{3}, \cite{4}, \cite{5}, \cite{6}, \cite{7}, and \cite{8}.
	Both the matrices, namely the distance Laplacian matrix and distance signless Laplacian matrix of a graph are positive semi-definite matrices. Consequently, both the matrices have non-negative eigenvalues.

	Over the last few decades, various researchers have pondered whether it is possible to predict the eigenvalues of a graph by observing the structure of a graph. One way to study the given problem is to perform various graph operations and create new graphs from existing graphs. Several graph operations have been introduced by researchers till now, some of them being \textit{join} of two graphs, \textit{disjoint union},
	\textit{Cartesian product}, \textit{direct product}, \textit{lexicographic product}.
	Several variants of corona product of two graphs have also been introduced and studied by various researchers in the recent past. Readers may refer to the papers \cite{9}, \cite{10}, \cite{11}, \cite{12}, \cite{13}, and \cite{14} for a detailed discussion in this regard.
	Moreover, researchers have determined the eigenvalues of the resulting graph operations in terms of existing graphs. Readers are suggested to see the papers \cite{15} and \cite{16} for more details.
	Recently, in \cite{17}, the authors have determined the  distance  Laplacian and distance signless Laplacian spectrum
	of \textit{generalized wheel graphs}.
	They have also introduced a new graph class and named it the  \textit{dumbbell graph.}
	The authors continued their study on dumbbell graphs in \cite{18}.
	The above works motivate us to study the Laplacian as well as the distance Laplacian spectrum of the \textit{generalized fan graph} in this paper. We have also introduced a new graph class and deduced its Laplacian and the distance Laplacian spectrum.

	\section{Preliminaries}
	\label{S2}
	
	The following definitions and theorems will be used in the subsequent sections. 
	
	\begin{definition}\cite{19}
		\label{EqP}
		Let $M$ be a order $n$ matrix defined as follows:
		\begin{center}
			\(
			\begin{pmatrix}
			M_{11}         & \cdots & M_{1t}     \\
			\vdots   & \ddots & \vdots \\
			M_{t1}            & \cdots & M_{tt}
			\end{pmatrix}.
			\)
		\end{center}
		Each block $M_{ij}$ has order  $n_i\times n_j$ for  $1\leq i, j\leq t$, and $M$ is equal to its transpose.
		Moreover, $n=n_1+\cdots+n_t$. For $1\leq i, j\leq t$, let $b_{ij}$ denote a matrix in which each element of $b_{ij}$ is obtained by adding all the entries in $M_{ij}$ and then dividing by the number of rows. The matrix  $B=(b_{ij})$ so obtained is known as the \textit{quotient} matrix of $M$. Additionally, if for each pair $i,j$, the sum of the entries in each row of $M_{ij}$ is constant,  then we call $B$ as the \textit{equitable quotient} matrix of $M$.
	\end{definition}
	There exists a relation between the set of eigenvalues of $B$ and $M$, which is given by the following theorem.
	\begin{theorem}\cite[Lemma  $2.3.1$]{19}
		\label{P1}
		If $\rho(M)$ is  the set of eigenvalues of $M$, and $\rho(B)$ is the set of eigenvalues of $B$, then $\rho(B)$ is contained in  $\rho(M)$. 
	\end{theorem}

	\section{Laplacian Spectra of Generalized Fan Graph and a New Graph Class}
	
	We first determine the eigenvalues of Laplacian matrix of 
	generalized fan graphs. We then introduce a new graph class and determine its Laplacian spectrum.

	\begin{definition}
		The generalized fan graph, denoted by $F_{m,n}$, is given by $F_{m,n}=\overline K_m+P_n$, where $\overline{K}_m$ is the null graph on $m$ vertices, and $P_n$ is the path graph on $n$ vertices.
	\end{definition}
	
	To determine the Laplacian spectrum of the generalized fan graph $F_{m,n}$, we shall first require the following result from \cite[Corollary 3.7]{20}.
	
	\begin{theorem}
		\label{Thjoin}
		Let $G_1+ G_2$ denote the join of two graphs $G_1$ and $G_2$.
		Then \begin{flalign*}
		\mu(G_1+ G_2;x)=\frac{x(x-n_1-n_2)}{(x-n_1)(x-n_2)}\mu(G_1,x-n_2)\mu(G_2,x-n_1),
		\end{flalign*}
		where $n_1$ and $n_2$ are orders of $G_1$ and $G_2$  respectively.
	\end{theorem}
	
	\begin{theorem}
		\label{II}
		If  $m,n\ge 2$, then the Laplacian eigenvalues  of $F_{m,n}$ are $0$ having multiplicity $1$, $m+n$ having multiplicity $1$, $n$ having multiplicity $m-1$, and $m+2-2\cos \frac{\pi j}{n}$  having multiplicity $1$ for $1\le j\le n-1$.
	\end{theorem}
	
	\begin{proof}
		
		We know that the Laplacian eigenvalues of $\overline K_m$ are $0$ having multiplicity $m$. 
		Hence,  $\mu(\overline{K}_m;x)=x^m$.
		Moreover, using \cite[Section 1.4.4]{19}, we find that the Laplacian eigenvalues of $P_n$ are $2-2\cos (\frac{\pi j}{n})$, where $ 0\le j\le n-1$.
		Hence, the characteristic polynomial of the Laplacian matrix of ${P}_n$ is given as follows:
		\begin{flalign*}
		\mu(P_n;x)&=x \times  \bigg[ \prod_{j=1}^{n-1}\bigg(x-2+2\cos \frac{\pi j}{n}\bigg)\bigg].
		\end{flalign*}
		
		Thus, using \Cref{Thjoin}, we get,
		\begin{flalign*}
		\mu(F_{m,n};x)&=\frac{x(x-m-n)}{(x-m)(x-n)}\times \mu(\overline{K}_m,x-n)\times \mu(P_n,x-m)
		\\
		&=\frac{x(x-m-n)}{(x-m)(x-n)}\times  (x-n)^m \times (x-m) \times  \bigg[ \prod_{j=1}^{n-1}\bigg(x-m-2+2\cos \frac{\pi j}{n}\bigg)\bigg]
		\\
		&=x(x-m-n)\times  (x-n)^{m-1}  \times  \bigg[ \prod_{j=1}^{n-1}\bigg(x-m-2+2\cos \frac{\pi j}{n}\bigg)\bigg].
		\end{flalign*}
		Hence the result follows.
	\end{proof}

	\begin{corollary}
		The Laplacian spectrum of the usual fan graph $F_{1,n}$ consists of $0$ having multiplicity $1$, $1+n$ having multiplicity $1$, and $3-2\cos \frac{\pi j}{n}$  having multiplicity $1$ for $1\le j\le n-1$.
	\end{corollary}

	\begin{proof}
		The proof follows from \cref{II} by putting $m=1$.
	\end{proof}

	We shall now introduce a new graph class and derive the Laplacian spectrum of the same. We shall denote the new graph class by $\mathcal{NC}(F_{m,n})$.
	We shall define the new graph in what follows.
	
	\begin{definition}
		\label{Def1}
		The graph $\mathcal{NC}(F_{m,n})$ has $2(m + n)$ vertices and is obtained by connecting $m$ vertices at the centers of two generalized fan graphs $F_{m,n}$, where $m,n \ge 2$ through $m$-edges. 
	\end{definition}
	
	We shall now illustrate the newly defined graph class $\mathcal{NC}(F_{m,n})$ with an example in what follows.
	
	\begin{example}
		We consider $m=3$ and $n=4$. We have the following two graphs namely, $\overline K_3$ and $P_3$.
		We shall first construct the generalized fan graph $F_{m,n}$.
		\begin{multicols}{2}
			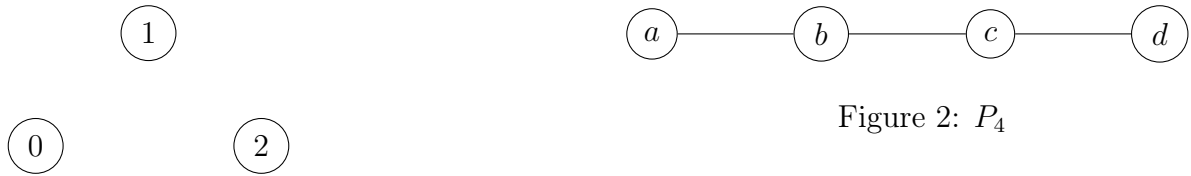
\begin{figure}[H]
				
				\begin{tikzpicture}[scale=0.5]

				\node[shape=circle,draw=black] (0) at (0,0) {$0$};  
				\node[shape=circle,draw=black] (1) at (3,3)  {$1$}; 
				\node[shape=circle,draw=black] (2) at (6,0)  {$2$}; 
				
				\end{tikzpicture}
				
				\caption{$\overline K_3$}
				\label{Figure 1}
			\end{figure}
			\begin{figure}[H]
				
				\begin{tikzpicture}[scale=0.75]

				\node[shape=circle,draw=black] (0) at (3,0) {$a$};  
				\node[shape=circle,draw=black] (1) at (6,0)  {$b$}; 
				\node[shape=circle,draw=black] (2) at (9,0)  {$c$}; 
				\node[shape=circle,draw=black] (3) at (12,0)  {$d$}; 
				\draw (0) -- (1); 
				\draw (1) -- (2); 
				
				\draw (2) -- (3);

				\end{tikzpicture}
				
				\caption{$P_4$}
				\label{Figure 2}
				
			\end{figure}
			
		\end{multicols}	 
		Using $\overline{K}_3$ and $P_4$, the generalized fan graph $F_{3,4}$ is given as follows:
		
		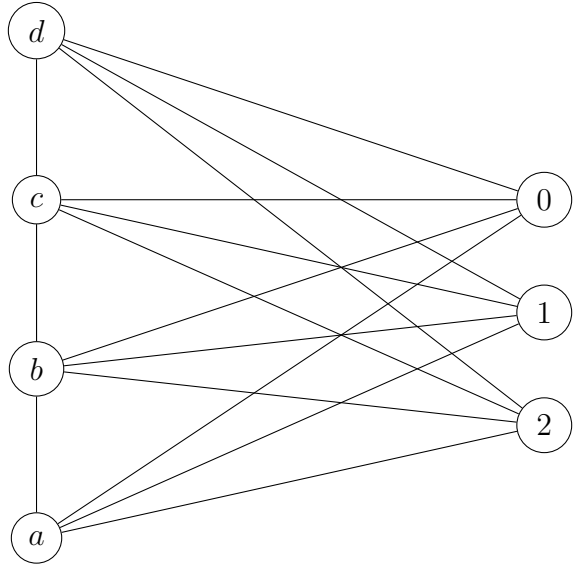
\begin{figure}[H]
			\centering
			\begin{tikzpicture}[scale=0.75]
			
			\node[shape=circle,draw=black] (0) at (0,3) {$a$};  
			\node[shape=circle,draw=black] (1) at (0,6)  {$b$}; 
			\node[shape=circle,draw=black] (2) at (0,9)  {$c$}; 
			\node[shape=circle,draw=black] (3) at (0,12)  {$d$}; 
			\node[shape=circle,draw=black] (a) at (9,9) {$0$};  
			\node[shape=circle,draw=black] (b) at (9,5)  {$2$}; 
			\node[shape=circle,draw=black] (c) at (9,7)  {$1$}; 
			\draw (0) -- (a); 
			\draw (0) -- (b); 
			\draw (0) -- (c);
			\draw (0) -- (1);
			\draw (1) -- (2);
			\draw (1) -- (2);
			\draw (2) -- (3); 
			\draw (1) -- (a); 
			\draw (1) -- (b); 
			\draw (1) -- (c);
			
			\draw (2) -- (a); 
			\draw (2) -- (b); 
			\draw (2) -- (c);
			
			\draw (3) -- (a); 
			\draw (3) -- (b); 
			\draw (3) -- (c);
			\end{tikzpicture}
			
			\caption{The generalized fan graph $F_{3,4}$.}
			\label{Figure 3}
			
		\end{figure}
		Using \Cref{Def1}, the new graph class $\mathcal{NC}(F_{3,4})$ is given as follows:
		
		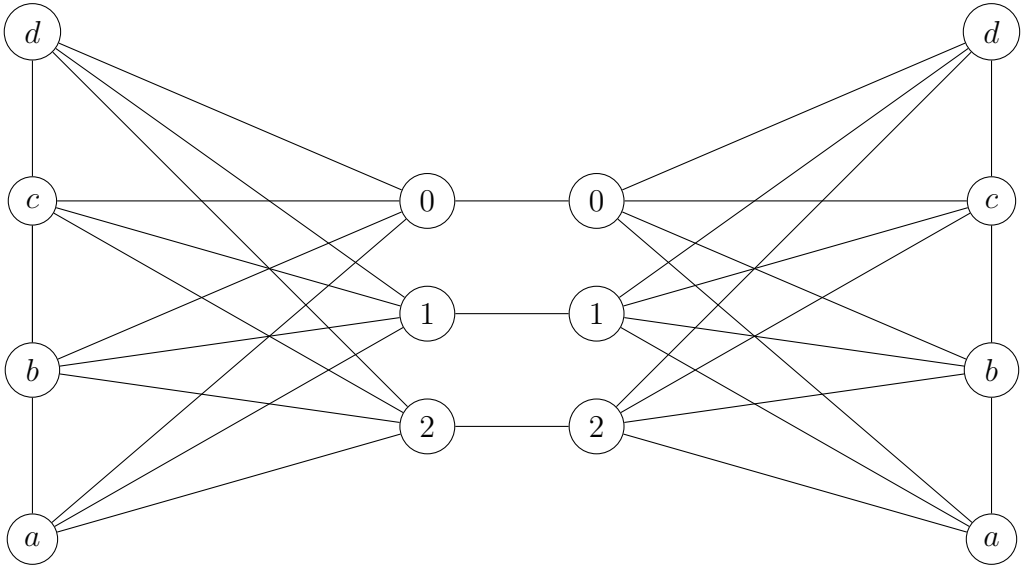
\begin{figure}[H]
			
			\begin{multicols}{2}
				
				\begin{tikzpicture}[scale=0.75]
				
				\node[shape=circle,draw=black] (0) at (2,3) {$a$};  
				\node[shape=circle,draw=black] (1) at (2,6)  {$b$}; 
				\node[shape=circle,draw=black] (2) at (2,9)  {$c$}; 
				\node[shape=circle,draw=black] (3) at (2,12)  {$d$}; 
				\node[shape=circle,draw=black] (a) at (9,9) {$0$};  
				\node[shape=circle,draw=black] (b) at (9,5)  {$2$}; 
				\node[shape=circle,draw=black] (c) at (9,7)  {$1$};
				\draw (0) -- (a); 
				\draw (0) -- (b); 
				\draw (0) -- (c);
				\draw (0) -- (1);
				\draw (1) -- (2);
				\draw (1) -- (2);
				\draw (2) -- (3); 
				\draw (1) -- (a); 
				\draw (1) -- (b); 
				\draw (1) -- (c);
				
				\draw (2) -- (a); 
				\draw (2) -- (b); 
				\draw (2) -- (c);
				
				\draw (3) -- (a); 
				\draw (3) -- (b); 
				\draw (3) -- (c);
				
				\node[shape=circle,draw=black] (a1) at (12,9) {$0$};  
				\node[shape=circle,draw=black] (b1) at (12,5)  {$2$}; 
				\node[shape=circle,draw=black] (c1) at (12,7)  {$1$};
				
				\node[shape=circle,draw=black] (01) at (19,3) {$a$};  
				\node[shape=circle,draw=black] (11) at (19,6)  {$b$}; 
				\node[shape=circle,draw=black] (21) at (19,9)  {$c$}; 
				\node[shape=circle,draw=black] (31) at (19,12)  {$d$};

				\draw (01) -- (a1); 
				\draw (01) -- (b1); 
				\draw (01) -- (c1);
				\draw (01) -- (11);
				\draw (11) -- (21);
				\draw (11) -- (21);
				\draw (21) -- (31); 
				\draw (11) -- (a1); 
				\draw (11) -- (b1); 
				\draw (11) -- (c1);
				
				\draw (21) -- (a1); 
				\draw (21) -- (b1); 
				\draw (21) -- (c1);
				
				\draw (31) -- (a1); 
				\draw (31) -- (b1); 
				\draw (31) -- (c1);
				
				\draw (a) -- (a1); 
				\draw (b) -- (b1); 
				\draw (c) -- (c1);
				
				\end{tikzpicture}
			\end{multicols}
			
			\caption{The graph $\mathcal{NC}_{3,4}$.}
			\label{Figure3}
		\end{figure}
		
	\end{example}
	We shall now illustrate the Laplacian eigenvalues of $\mathcal{NC}_{m,n}$ in what follows.
	It is known that the Laplacian eigenvalues of $P_n$ are $0$ and $2(1-\cos \frac{\pi j}{n})$ having multiplicity $1$ for $1\le j\le n-1$.

	\begin{theorem}
		\label{I}
		If $m,n\ge 2$, then the Laplacian eigenvalues of $\mathcal{NC}(F_{m,n})$ are 
		as follows:
		\begin{enumerate}
			\item [$\bullet$] $2(1-\cos \frac{\pi j}{n})+m$ having multiplicity $2$ for $1\le j\le n-1$,
			\item [$\bullet$] $n$ having multiplicity $m-1$,
			\item [$\bullet$] $n+2$ having multiplicity $m-1$,
			\item [$\bullet$] $\frac{m+n}{2} \pm \frac{\sqrt{(m^2 + 2(m + 2)n + n^2 - 4m + 4) + 1}}{2}$ having multiplicity $1$,
			\item [$\bullet$]$m+n$ having multiplicity $1$,
			\item [$\bullet$] $0$ having multiplicity $1$.
		\end{enumerate}
		
	\end{theorem}
	
	\begin{proof}
			
		We shall first index the vertices of $P_n$, then list the vertices of $\overline{K}_m$. We again list the vertices of the second copy of  $\overline{K}_m$ and finally list the vertices of the second copy of $P_n$.
		Thus the Laplacian matrix of $\mathcal{NC}(F_{m,n})$ is given as follows:
		\begin{flalign*}
		L(\mathcal{NC}(F_{m,n}))=
		\left(\begin{matrix}
		L(P_n)+mI && -J_{n\times m} && 0_{n\times m} && 0_{n\times n}
		\\
		\\
		-J_{m\times n} && (n+1)I_{m\times m} && -I_{m\times m} && 0_{m\times n}
		\\
		\\
		0_{n\times m} && -I_{m\times m} && (n+1)I_{m\times m} && -J_{m\times n}
		\\
		\\
		0_{n\times n}&&  0_{n\times m} && -J_{n\times m} && L(P_n)+mI
		\end{matrix}\right).
		\end{flalign*}

		Now, since $L(P_n)$ is a singular matrix, so zero will be an eigenvalue of $L(P_n)$.
		The eigenvector corresponding to the eigenvalue $0$ is $\mathbf{1}=[1,1,\dots, 1]^T$.
		For a symmetric matrix, if $\lambda_i$ and $\lambda_j$ are two distinct eigenvalues with eigenvectors $v_i$ and $v_j$ respectively, then $v_i$ and $v_j$ are orthogonal to each other.
		Let $\lambda(\neq 0)$ be an eigenvalue of $L(P_n)$ having  eigenvector $\mathbf{v}$.
		Then, $\mathbf{1}^T\mathbf{v}=0$.

		Let $v_i$, $2\le i\le m$ be an eigenvector corresponding to the eigenvalue $\lambda_i=2(1-\cos \frac{\pi i}{n})$ of $P_n$.
		Let $\mathbf{V_i}=\left(\begin{array}{cc}
		\mathbf{v_i}_{n}\\ \mathbf{0}_{m}\\ \mathbf{0}_{m}\\\mathbf{0}_{n}
		\end{array}\right)$.
		Now	$L(\mathcal{NC}(F_{m,n}))\mathbf{V_i}= (\lambda_i+m)\mathbf{V_i}$.
		Thus,  $\lambda_i+m$ is an eigenvalue of $L(\mathcal{NC}(F_{m,n}))$.
		Similarly, let $\mathbf{W_i}=\left(\begin{array}{cc}
		\mathbf{0}_{n}\\ \mathbf{0}_{m}\\ \mathbf{0}_{m}\\\mathbf{v_i}_{n}
		\end{array}\right)$, we observe that  	$L(\mathcal{NC}(F_{m,n}))\mathbf{W_i}= (\lambda_i+m)\mathbf{W_i}$.
		Thus, again, we find that $\lambda_i+m$ is an eigenvalue of $L(\mathcal{NC}(F_{m,n}))$ for $2\le i\le m$.
		Hence, we observe that $\lambda_i+m$ is an eigenvalue of $L(\mathcal{NC}(F_{m,n}))$ for $2\le i\le m$ having multiplicity $2$.
		
		Let $\mathbf{X_i}=\left(\begin{array}{cc}
		\mathbf{0}_{n}\\ \mathbf{v_i}_{m}\\ \mathbf{v_i}_{m}\\\mathbf{0}_{n}
		\end{array}\right)$.
		
		We have
		\begin{flalign*}
		&L(\mathcal{NC}(F_{m,n}))\mathbf{X_i}
		\\
		&=\left(\begin{matrix}
		L(P_n)+mI && -J_{n\times m} && 0_{n\times m} && 0_{n\times n}
		\\
		\\
		-J_{m\times n} && (n+1)I_{m\times m} && -I_{m\times m} && 0_{m\times n}
		\\
		\\
		0_{n\times m} && -I_{m\times m} && (n+1)I_{m\times m} && -J_{m\times n}
		\\
		\\
		0_{n\times n}&&  0_{n\times m} && -J_{n\times m} && L(P_n)+mI
		\end{matrix}\right)
		\left(\begin{array}{cc}
		\mathbf{0}_{n}\\\\ \mathbf{v_i}_{m}\\\\ \mathbf{v_i}_{m}\\\\\mathbf{0}_{n}
		\end{array}\right)
		\\
		&=\left(\begin{array}{cc}
		\mathbf{0}\\\\((n+1)-1)\mathbf{v_i}_{m}\\\\ ((n+1)-1)\mathbf{v_i}_{m}\\\\\mathbf{0}
		\end{array}\right)
		\\
		&=\left(\begin{array}{cc}
		\mathbf{0}\\\\n\mathbf{v_i}_m\\\\ n\mathbf{v_i}_m\\\\\mathbf{0}
		\end{array}\right)
		\\
		&=n\left(\begin{array}{cc}
		\mathbf{0}\\\\\mathbf{v_i}_{m}\\\\ \mathbf{v_i}_{m}\\\\\mathbf{0}
		\end{array}\right).
		\end{flalign*}
		
		We thus obtain  $L(\mathcal{NC}(F_{m,n}))\mathbf{X_i}= n\mathbf{X_i}$.
		Thus,  $n$ is an eigenvalue of $L(\mathcal{NC}(F_{m,n}))$.
		Hence, we find that $n$ is an eigenvalue of $L(\mathcal{NC}(F_{m,n}))$  having multiplicity $m-1$.
		
		Let $\mathbf{Y_i}=\left(\begin{array}{cc}
		\mathbf{0}_{n}\\ \mathbf{v_i}_{m}\\ \mathbf{-v_i}_{m}\\\mathbf{0}_{n}
		\end{array}\right)$.
		Now $L(\mathcal{NC}(F_{m,n}))\mathbf{X_i}= (n+2)\mathbf{Y_i}$.
		Thus,  $n+2$ is an eigenvalue of $L(\mathcal{NC}(F_{m,n}))$ having multiplicity $m-1$.
		
		Thus, we determine $2(n+m-2)$ eigenvalues of $L(\mathcal{NC}(F_{m,n})$.
		We shall now use \Cref{EqP}.
		We shall now use \Cref{P1} to find the  $4$ remaining  eigenvalues of $L(\mathcal{NC}(F_{m,n})$. We find that they are contained in the spectrum of matrix $B$ given as follows:
		
		\[
		B=
		\left(
		\begin{array}{cccccccc}
		m &&-m && 0 && 0
		\\
		\\
		-n && n+1 && -1 && 0
		\\
		\\
		0 && -1 && n+1 && -n
		\\
		\\
		0 && 0 && -m && m 
		\end{array}
		\right).
		\]
		The characteristic polynomial  of $B$ is :
		\begin{flalign*}
		\Theta(B,x)&=x^4 + (-2m - 2n - 2)x^3 + (m^2 + 2mn + n^2 + 4m + 2n)x^2 + (-2m^2 - 2mn)x.
		\end{flalign*}
		
		On solving $\Theta(B,x)=0$, we obtain the required result.

	\end{proof}

	\section{Distance Laplacian Spectrum of Generalized Fan Graph and a New Graph Class}
	\label{S3}
	In this section, we evaluate the distance Laplacian spectrum of 
	the generalized fan graph. We then determine the distance Laplacian spectrum of the new graph class that was introduced in the previous section.	
	To determine the distance Laplacian spectrum of the
	generalized fan graph, we shall need the given theorem.

	\begin{theorem}\label{Th1}
		\label{Join}
		Let $G_1$ be a graph on $n_1$ vertices  having Laplacian eigenvalues $0=\lambda_1\le \lambda_2\le\cdots \le \lambda_{n_1}$ and  $G_2$ be a graph on $n_2$ vertices  having Laplacian eigenvalues $0=\mu_1\le \mu_2\le\cdots \le \mu_{n_2}$. Then the distance Laplacian spectrum of $G_1+ G_2$ consists of  $n_2+2n_1-\lambda_i$ having multiplicity $1$ for $2\le i\le n_1$,
		$n_1+2n_2-\mu_i$ having multiplicity $1$ for $2\le j\le n_2$, and $0, n_1+n_2$ having multiplicity $1$.
		
	\end{theorem}
	
	\begin{proof}
		We shall first index the vertices of the graph $G_1$.
		We then index the vertices of the graph $G_2$.
		We have:
		\begin{flalign*}
		D^L(G_1+ G_2)&=
		\left(\begin{matrix}
		D^{L_1} && -J_{n_1\times n_2} 
		\\
		\\
		-J_{n_1\times n_2} && D^{L_2}
		\end{matrix}\right).
		\end{flalign*}
		Here,
		\begin{flalign*}
		D^{L_1}&=Tr(G_1)-D(G_1)
		\\
		&=Tr(G_1)+A(G_1)-2J_{n_1\times n_1}+2I_{n_1\times n_1}
		\\
		&=\bigg((n_2+2(n_1-1))I_{n_1\times n_1}\bigg)-\text{Deg}(G_1)
		\\&+A(G_1)-2J_{n_1\times n_1}+2I_{n_1\times n_1}
		\\
		&=\bigg((n_2+2(n_1-1)+2)I_{n_1\times n_1}\bigg)-\text{Deg}(G_1)+A(G_1)-2J_{n_1\times n_1}
		\\
		&=\bigg((n_2+2n_1)I_{n_1\times n_1}\bigg)-\text{Deg}(G_1)+A(G_1)-2J_{n_1\times n_1}
		\\
		&=\bigg((n_2+2n_1)I_{n_1\times n_1}\bigg)-2J_{n_1\times n_1}-L(G_1),
		\end{flalign*}
		and,	
		\begin{flalign*}
		D^{L_2}&=Tr(G_2)-D(G_2)
		\\
		&=\bigg((n_1+2n_2)I_{n_2\times n_2}\bigg)-2J_{n_2\times n_2}-L(G_2).
		\end{flalign*}
		
		We know that the Laplacian matrix  $L(G_1)$ is a singular matrix having a determinant as $0$.
		Moreover, since the sum of the entries of each row is $0$, so $0$ will be an eigenvalue of $L(G_1)$.
		Hence, we have $L(G_1)\mathbf{1}=L(G_1)[1,1,\dots, 1]^T=\mathbf{0}$.
		Let $\lambda_i$ be a non-zero eigenvalue of $L(G_1)$ whose eigenvector is $\mathbf{v_i}$, $2\le i\le n$. 	Moreover, $\mathbf{1}^T\mathbf{v_i}=0$.

		Let  $\mathbf{V_i}=\left(\begin{array}{cc}
		\mathbf{v_i}_{n_1}\\ \mathbf{0}_{n_2}
		\end{array}\right)$.
		We obtain, 
		\begin{flalign*}
		&D^L(G_1+ G_2)\mathbf{V_i}
		\\
		&=\left(\begin{matrix}
		D^{L_1} & -J_{n_1\times n_2} \\
		\\
		-J_{n_2\times n_1} & D^{L_2}
		\end{matrix}\right)\left(\begin{array}{cc}
		\mathbf{v_i}_{n_1}\\\\ \mathbf{0}_{n_2}
		\end{array}\right)
		\\
		&=\left(\begin{array}{cc}
		D^{L_1}\mathbf{v_i}\\\\ \mathbf{0}
		\end{array}\right)
		\\
		&=\left(\begin{array}{cc}\bigg(((n_2+2n_1)I_{n_1\times n_1})-2J_{n_1\times n_1}-L(G_1)\bigg)\mathbf{v_i}\\\\ \mathbf{0}\end{array}\right)
		\\		
		&=\left(\begin{array}{cc}(n_2+2n_1)\mathbf{v_i}-\lambda_i\mathbf{v_i}\\\\ \mathbf{0}\end{array}\right)
		\\
		&=\left(\begin{array}{cc}(n_2+2n_1-\lambda_i)\mathbf{v_i}\\\\ \mathbf{0}\end{array}\right)
		\\
		&=(n_2+2n_1-\lambda_i)\mathbf{V_i}.
		\end{flalign*}
		
		Thus,  if  $\lambda_i$ is an eigenvalue  of $L(G_1)$ for $2\le i\le n_1$, we find  that  $n_2+2n_1-\lambda_i$ is an eigenvalue of $D^L(G_1+ G_2)$.
		This provides us with $n_1-1$ distance Laplacian eigenvalues of  $G_1+G_2$.
		
		Let $\mu_j$ be an eigenvalue  of $L(G_2)$.
		Let  $\mathbf{w}$ be an eigenvector of $\mu_j$.
		Using similar arguments as given above, we find that  $n_1+2n_2-\mu_i$ is a distance Laplacian eigenvalue of $G_1+ G_2$ corresponding to eigenvector $\mathbf{W}$.
		Here, $\mathbf{W}=\left(\begin{array}{cccccccc}
		\mathbf{0}_{n_1}\\\mathbf{w}_{n_2}
		\end{array}\right).$
		This provides us with $n_1+n_2-2$ distance Laplacian eigenvalues of $G_1+G_2$.
		The remaining two eigenvalues of $D^L(G_1+G_2)$ can be obtained by using the concept of equitable partitions(\Cref{EqP}).
		Since each block matrix of $D^L(G_1+ G_2)$ has a constant row sum, we find that the equitable quotient matrix of $D^L(G_1+ G_2)$ is given as follows:
		\begin{flalign*}
		B&=\left(
		\begin{array}{cccc}
		n_2&& -n_2\\
		-n_1&&n_1
		\end{array}
		\right).
		\end{flalign*}
		Since $\sigma(B)=\left(\begin{array}{ccccc}
		n_1+n_2 & & 0\\
		1 && 1
		\end{array}\right)$, using Theorem \ref{P1}, we find that the eigenvalues of $D^L(G_1+ G_2)$ are $n_2+2n_1-\lambda_i$ having multiplicity $1$ for $2\le i\le n_1$,
		$n_1+2n_2-\mu_i$ having multiplicity $1$ for $2\le j\le n_2$, and $0, n_1+n_2$ having multiplicity $1$.
		
	\end{proof}
	
	We now determine the distance Laplacian spectra of the generalized fan graph $F_{m,n}$.
	
	\begin{theorem}
		\label{Fan1}
		The  spectrum of the distance Laplacian matrix of  $F_{m,n}$ consists of $n+m$ having multiplicity $m-1$, $m+2n-2+2\cos (\frac{\pi j}{n})$ having multiplicity $1$ for $0\le j\le n-1$, and $0,m+n$ having multiplicity $1$.
	\end{theorem}
	
	\begin{proof}
		We know $F_{m,n}=\overline K_m+P_n$.
		Using \Cref{Th1}, the  eigenvalues of the distance Laplacian matrix of  $F_{m,n}$ are $n+m$ having multiplicity $m-1$, $m+2n-2+2\cos (\frac{\pi j}{n})$ having multiplicity $1$ for $0\le j\le n-1$, and $0,m+n$ having multiplicity $1$.
		
	\end{proof}
	
	\begin{corollary}
		The distance  Laplacian spectrum  of the usual fan graph $F_{1,n}$ consists of  $2n-1+2\cos (\frac{\pi j}{n})$ having multiplicity $1$ for $0\le j\le n-1$, and $0,n+1$ having multiplicity $1$.
	\end{corollary}
	
	\begin{proof}
		The proof follows by substituting $m=1$ in \Cref{Fan1}.
	\end{proof}

	\subsection{Distance Laplacian spectrum of $\mathcal{NC}(F_{m,n})$}

	In our next theorem, we shall now determine the distance Laplacian spectrum of the new graph class $\mathcal{NC}(F_{m,n})$.

	\begin{theorem}
		\label{1}
		The distance Laplacian spectrum  of  $\mathcal{NC}(F_{m,n})$ consists of:
		\begin{itemize}
			\item [$\bullet$]  $(5n+3m-\lambda_i)$ having multiplicity $2$ for $2\le i\le n$,
			\item [$\bullet$] $3n+5m-4$ having multiplicity $m-1$,
			\item [$\bullet$]  $3n+5m-2$ having multiplicity $m-1$,
			\item[$\bullet$]	$\frac{9}{2}(n +m) - 2\pm \frac{1}{2}\sqrt{A} $, where $A=24n + 9n^2 - 14nm-24m + 9m^2 + 16$, having multiplicity $1$,
			\item[$\bullet$] $3(n+m)$ having multiplicity $1$,
			and \item[$\bullet$] $0$ having multiplicity $1$.
		\end{itemize}
		
	\end{theorem}
	
	\begin{proof}
		We shall first index the vertices of $P_n$, then list the vertices of $\overline{K}_m$. We again list the vertices of the second copy of  $\overline{K}_m$ and finally list the vertices of the second copy of $P_n$.
		We have:
		\begin{flalign*}
		&D^L(\mathcal{NC}(F_{m,n}))
		\\&=
		\left(\begin{matrix}
		D^{L_1} && -J_{n\times m} && -2J_{n\times m} && -3J_{n\times m}
		\\
		\\
		-J_{m\times n} && D^{L_2} && -(3J-2I)_{m\times m} && -2J_{m\times n}
		\\
		\\
		-2J_{m\times n} && -(3J-2I)_{m\times m} && D^{L_2} && -J_{m\times n}
		\\
		\\
		-3J_{n\times n}&&  -2J_{n\times m} && -J_{n\times m} && D^{L_1}
		\end{matrix}\right).
		\end{flalign*}
		Here,
		
		\begin{flalign*}
		D^{L_1}&=(5n+3m)I_{n}-2J_{n\times n}-L(G_1), 	\text{ and }
		\\
		D^{L_2}&=(3n+5m-2)I_{m}-2J_{m\times m}.
		\end{flalign*}

		Assuming $\lambda_i$ to be an eigenvalue of $L(P_n)$ with eigenvector $\mathbf{v_i}$ for  $2\le i\le n$, we have $\mathbf{1}^T\mathbf{v_i}=0$.
		
		Considering $\mathbf{V_i}=\left(\begin{array}{cc}
		\mathbf{v_i}_{n}\\ \mathbf{0}_{m}\\ \mathbf{0}_{m}\\\mathbf{0}_{n}
		\end{array}\right)$, we obtain, 
		\begin{flalign*}
		&D^L(\mathcal{NC}(F_{m,n}))
		\\
		&=\left(\begin{matrix}
		D^{L_1} & -J_{n\times m} & -2J_{n\times m} & -3J_{n\times n}
		\\
		\\
		-J_{m\times n} & D^{L_2} & -(3J-2I)_{m\times m} & -2J_{m\times n}
		\\
		\\
		-2J_{m\times n} & -(3J-2I)_{m\times m} & D^{L_2} & -J_{m\times n}
		\\
		\\
		-3J_{n\times m}&  -2J_{n\times m} & -J_{n\times m} & D^{L_1}
		\end{matrix}\right)\left(\begin{array}{cc}
		\mathbf{v_i}_{n}\\\\ \mathbf{0}_{m}\\\\ \mathbf{0}_{m}\\\\\mathbf{0}_{n}
		\end{array}\right)
		\\
		&=\left(\begin{array}{cc}
		D^{L_1}\mathbf{v_i}\\\\ \mathbf{0}\\\\ \mathbf{0}\\\\\mathbf{0}
		\end{array}\right)
		\\
		&=\left(\begin{array}{cc}\bigg((5n+3m)I_{n}-2J_{n\times n}-L(P_n)\bigg)\mathbf{v_i}\\\\ \mathbf{0}\\\\ \mathbf{0}\\\\\mathbf{0}\end{array}\right)
		\\		
		&=\left(\begin{array}{cc}(5n+3m)\mathbf{v_i}-L(P_n)\mathbf{v_i}\\\\ \mathbf{0}\\\\ \mathbf{0}\\\\\mathbf{0}\end{array}\right)
		\\
		&=(5n+3m-\lambda_i)\mathbf{V_i}.
		\end{flalign*}
		
		Thus, we observe that $(5n+3m-\lambda_i)$ becomes an eigenvalue of $D^L(\mathcal{NC}(F_{m,n}))$ having multiplicity $1$. Here, $2\le i\le n$.
		
		Let  $\mathbf{W_i}=\left(\begin{array}{cc}
		\mathbf{0}_{n}\\ \mathbf{0}_{m}\\ \mathbf{0}_{m}\\\mathbf{v_i}_{n}
		\end{array}\right)$ for $2\le i\le n$. We observe that $D^L(\mathcal{NC}(F_{m,n}))\mathbf{W_i}=(5n+3m-\lambda_i)\mathbf{W_i}.$
		Thus, $(5n+3m-\lambda_i)$ is an eigenvalue of $D^L(\mathcal{NC}(F_{m,n}))$ having multiplicity $1$ for $2\le i\le n$.
		Hence we find that $(5n+3m-\lambda_i)$ is an eigenvalue of $D^L(\mathcal{NC}(F_{m,n}))$ having multiplicity $2$ for each $2\le i\le n$.
		
		Moreover, we observe that $3n+5m$ and $3n+5m-4$  are eigenvalues of $D^L(\mathcal{NC}(F_{m,n}))$ having multiplicity $m-1$.
		
		This gives us $2(n+m-2)$ eigenvalues of $D^L(\mathcal{NC}(F_{m,n}))$.
		The remaining $4$  eigenvalues of $D^L(\mathcal{NC}(F_{m,n}))$ are contained in the spectrum of $B$ where
		
		\[
		B=
		\left(
		\begin{array}{cccccccc}
		3(n+m) &&-m && -2m && -3n
		\\
		\\
		-n && 3(n+m)-2 && -(3m-2)&& -2n
		\\
		\\
		-2n && -(3m-2)&& 3(n+m)-2&& -n
		\\
		\\
		-3n&& -2m&& -m && 3(n+m) 
		\end{array}
		\right).
		\]
		The characteristic polynomial of $B$ is $x^4 + (-12m - 12n + 4)x^3 + (45m^2 + 98mn + 45n^2 - 24m - 36n)x^2 + (-54m^3 - 186m^2n - 186mn^2 - 54n^3 + 36m^2 + 108mn + 72n^2)x$.
		
		On solving, we find that the eigenvalues of $B$ are	$\frac{9}{2}(n +m) - 2\pm \frac{1}{2}\sqrt{A}$, $3(n + m)$ and $0$ having multiplicity $1$, where $A=24n + 9n^2 - 14nm-24m + 9m^2 + 16$.

	\end{proof}

	\section{Conclusion}
	
	In this paper, our main aim is to determine the Laplacian and the distance Laplacian spectrum of the generalized fan graph $F_{m,n}$.
	Moreover, we also introduce a new graph class
	and	determine its Laplacian as well as its distance Laplacian spectrum.
	We illustrate our results with various examples.
	
	We now provide a problem to the readers for future work.
	The matrix $D_t(G) = t\text{Tr}(G) + (1-t) D(G), 0<t<1$, is known as the  \textit{generalized distance matrix} of a graph $G$.
	We encourage the readers to determine the spectrum of the generalized distance matrix of the generalized fan graph as well as the new graph class introduced in this paper.
	
	\section{Declarations}
	
	\subsection{Conflict of interest:}  The authors state that there is no conflict of interest.
	
	\subsection{Funding:}   Not Applicable.
	
	\subsection{Authors contribution:}   
	S.B.-Conceptualization, Methodology, Formal Analysis, Writing-Review and Editing, Supervision. 
	S.G.- Resources, Writing-Original Draft, Methodology.
	The final submitted	version of this manuscript has been read and approved by all the authors.

	\subsection{Acknowledgement:} Not Applicable.
	
	\subsection{Data availability statement:} The article includes all the data that support the findings of this study.

	\section{Appendix}
	
	In this section, we provide the adjacency and the Laplacian eigenvalues of Fan graph $F_{1,n}$(\Cref{Tab}) and the generalized Fan graph $F_{m,n}$(\Cref{Tab1}) for various values of $m$ and $n$.
	
		\begin{table}[ht]
		\begin{center}
			\begin{tabular}{c|c|c}
				$n$ & $\text{Adjacency Spectrum}$ & $\text{Laplacian Spectrum}$ \\
				\hline
				$3$ &   $\{0, -1, -1.56, 2.56\}$ & $\{0, 1, 1,4\}$  \\
				\hline
				$4$ & $\{-1.62, 0.62, -1.47, -0.46, 2.93\}$ & $\{5, 3, 0, 1.58, 4.41\}$  \\
				\hline
				$5$ & $\{3.22, 0.11, -1.53, -1.81,  1,-1\}$ & $\{6, 0, 2.38, 4.62, 1.38, 3.62\}$ \\
				\hline
				$6$ & $\{-1.80, -0.44, 1.25, -1.82, -1.18, 0.54, 3.46\}$ & $\{7, 4, 3, 2, 0, 1.27, 4.73\}$\\
				\hline
				$7$ & $\{0, -2, -1.41, 1.41, -1.81, -0.71, 0.84, 3.67\}$ & $\{8, 0, 1.75, 3.44, 4.8, 1.2, 2.55, 4.25\}$\\
			\end{tabular}
		\end{center}
		
		\caption{Comparison Table for Adjacency Spectrum \& Laplacian Spectrum of Fan Graph }
		\label{Tab}
	\end{table}

		\begin{table}[ht]
		\begin{center}
			\begin{tabular}{c|c|c|c}
				$n$ & $m$ &$\text{Adjacency Spectrum}$ & $\text{Laplacian Spectrum}$ \\
				\hline
				$2$ & $2$&   $\{2, -2, 0, 0\}$ & $\{0, 2, 2,4\}$  \\
				\hline
				$2$ &$3$ &$\{-2, 0, 0, -1.24, 3.24\}$ & $\{0, 5, 5, 3, 3\}$  \\
				\hline
				$3$ &$2$ & $\{3, -1, -2, 0, 0\}$ & $\{0,5,5,2,2\}$ \\
				\hline
				$3$ & $4$ &$\{0, 0, -1.62, 0.62, -2.84, -0.49, 4.32\}$ & $\{0,7, 5, 4, 4, 3.58, 6.41\}$\\
				\hline
				$4$ &$3$ & $\{0, 0, 0, 0, -2.92, -1.3, 4.22\}$ & $\{0, 5, 7, 7, 3, 3, 3\}$\\
			\end{tabular}
		\end{center}
		
		\caption{Comparison Table for Adjacency Spectrum \& Laplacian Spectrum of Generalized Fan Graph }
		\label{Tab1}
	\end{table}

\end{document}